\documentclass[12pt]{article}

\usepackage{latexsym,amsmath}
\usepackage{amssymb}
\usepackage{url}
\usepackage{graphicx}
\usepackage{microtype}

\usepackage[margin=1.2in]{geometry}

\newtheorem{thm}{Theorem}[section]

\newtheorem{lemma}[thm]{Lemma}
\newtheorem{prop}[thm]{Proposition}

\newtheorem{conj}[thm]{Conjecture}

\newtheorem{question}[thm]{Question}

\newcommand{\qed}[0]{{\hspace*{\fill}\mbox{$\Box$}}}

\begin{document}

\renewcommand{\thefootnote}{\fnsymbol{footnote}}

\title{Counting colorings of a regular graph}

\author{David Galvin\thanks{Department of Mathematics,
University of Notre Dame, 255 Hurley Hall, Notre Dame IN
46556; dgalvin1@nd.edu.}}

\maketitle

\begin{abstract}
At most how many (proper) $q$-colorings does a regular graph admit? Galvin and Tetali conjectured that among all $n$-vertex, $d$-regular graphs with $2d|n$, none admits more $q$-colorings than the disjoint union of $n/2d$ copies of the complete bipartite graph $K_{d,d}$. In this note we give asymptotic evidence for this conjecture, showing that the number of proper $q$-colorings admitted by an $n$-vertex, $d$-regular graph is at most
$$
\begin{array}{ll}
\left(q^2/4\right)^\frac{n}{2} \binom{q}{q/2}^\frac{n(1+o(1))}{2d} & \mbox{if $q$ is even}, \\
\left((q^2-1)/4\right)^\frac{n}{2} \binom{q+1}{(q+1)/2}^\frac{n(1+o(1))}{2d} & \mbox{if $q$ is odd},
\end{array}
$$
where $o(1)\rightarrow 0$ as $d \rightarrow \infty$; these bounds agree up to the $o(1)$ terms with the counts of $q$-colorings of $n/2d$ copies of $K_{d,d}$.

An auxiliary result is an upper bound on the number of colorings of a regular graph in terms of its independence number. For example, we show that for all even $q$ and fixed $\varepsilon > 0$ there is $\delta=\delta(\varepsilon,q)$ such that the number of proper $q$-colorings admitted by an $n$-vertex, $d$-regular graph with no independent set of size $n(1-\varepsilon)/2$ is at most
$$
\left(q^2/4-\delta\right)^\frac{n}{2},
$$
with an analogous result for odd $q$.
\end{abstract}

\section{Introduction}

Throughout, $G$ is a simple, finite loopless graph, and $q$ is a positive integer. A {\em proper $q$-coloring} (or just {\em $q$-coloring}) of $G$ is a function from the vertices of $G$ to $\{1, \ldots, q\}$ with the property that adjacent vertices have different images. We write ${\mathcal c}_q(G)$ for the number of $q$-colorings of $G$.

The following is a natural extremal enumerative question: for a family ${\mathcal G}$ of graphs, which $G \in {\mathcal G}$ maximizes $c_q(G)$? For example, for the family of $n$-vertex, $m$-edge graphs this question was raised independently by Wilf \cite{BenderWilf,Wilf} and Linial \cite{Linial}, who both came across it in their study of running times of coloring algorithms. Although it has only been answered completely in some very special cases many partial results have been obtained (see \cite{LohPikhurkoSudakov} for a good history of the problem).

The focus of this note is the family ${\mathcal G}(n,d)$ of $n$-vertex $d$-regular graphs with $d \geq 2$ (the case $d=1$ being trivial). Galvin and Tetali \cite{GalvinTetali-weighted} used an entropy argument to show that for $2d|n$ no {\em bipartite} $G$ in ${\mathcal G}(n,d)$ admits more $q$-colorings, for each $q\geq 2$, than $\frac{n}{2d}K_{d,d}$, the disjoint union of $n/2d$ copies of the complete bipartite graph $K_{d,d}$ with $d$ vertices in each partite set. More generally they found $c_q(G) \leq c_q(K_{d,d})^{n/2d}$ for all $n$, $d$ and bipartite $G \in {\mathcal G}(n,d)$ (this is
\cite[Prop. 1.2]{GalvinTetali-weighted} in the special case $H=K_q$), and they conjectured that this bound should still hold when the biparticity assumption is dropped.
\begin{conj} \label{conj-qcols}
Fix $d \geq 2$ and $n \geq d+1$. For any $G \in {\mathcal G}(n,d)$ and any $q \geq 2$,
$$
c_q(G) \leq c_q(K_{d,d})^\frac{n}{2d}.
$$
\end{conj}
For $q=2$ this follows immediately from the bipartite case established in \cite{GalvinTetali-weighted}. Zhao \cite{Zhao2} established the conjecture for all $q \geq (2n)^{2n+2}$, and in the case $2d|n$ Galvin \cite{Galvin-homsconjcorr}, using ideas introduced by Lazebnik \cite{Lazebnik} on a related problem, reduced this to $q > 2\binom{nd/2}{4}$, but neither the approach of \cite{Galvin-homsconjcorr} nor that of \cite{Zhao2} seems adaptable to the case of constant $q \geq 3$.

Conjecture \ref{conj-qcols} is a special case of a more general conjecture concerning graph homomorphisms. A {\em homomorphism} from $G$ to a graph $H$ (which may have loops) is a map from vertices of $G$ to vertices of $H$ with adjacent vertices in $G$ being mapped to adjacent vertices in $H$. Homomorphisms generalize $q$-colorings (if $H=K_q$ then the set of homomorphisms to $H$ is in bijection with the set of $q$-colorings of $G$) as well as other graph theory notions, such as independent sets. A {\em independent set} in a graph is a set of pairwise non-adjacent vertices; notice that if $H=H_{\rm ind}$ is the graph on two adjacent vertices with a loop at exactly one of the vertices, then a homomorphism from $G$ to $H$ may be identified, via the preimage of the unlooped vertex, with an independent set in $G$. Amending a false conjecture from \cite{GalvinTetali-weighted}, the following conjecture is made in \cite{Galvin-homsconjcorr}. Here we write ${\rm hom}(G,H)$ for the number of homomorphisms from $G$ to $H$.
\begin{conj} \label{conj-homs}
Fix $d \geq 2$ and $n \geq d+1$. For any $G \in {\mathcal G}(n,d)$ and any finite graph $H$ (perhaps with loops, but without multiple edges),
$$
{\rm hom}(G,H) \leq \max \left\{{\rm hom}(K_{d,d},H)^\frac{n}{2d}, {\rm hom}(K_{d+1},H)^\frac{n}{d+1} \right\},
$$
where $K_{d+1}$ is the complete graph on $d+1$ vertices.
\end{conj}
When $d \geq q$ we have ${\rm hom}(K_{d+1},K_q)=0$ and so in this range Conjecture \ref{conj-homs} implies Conjecture \ref{conj-qcols}.

\medskip

The inspiration for Conjecture \ref{conj-homs}, and the partial result of \cite{GalvinTetali-weighted} that the conjecture is true for all {\em bipartite} $G$, was the special case of enumerating independent sets ($H=H_{\rm ind})$. In what follows we use $i(G)$ to denote the number of independent sets in $G$. Alon \cite{Alon} conjectured that for all $G \in {\mathcal G}(n,d)$ we have 
$$
i(G) \leq i(K_{d,d})^{n/2d} = (2^{d+1}-1)^{n/2d} = 2^{n/2 + O(n(1+o(1))/2d)},
$$ 
and proved the weaker bound $i(G) \leq 2^{n/2 + O(n/d^{1/10})}$. The sharp bound was proved for {\em bipartite} $G$ by Kahn \cite{Kahn}, but it was a while before a bound for general $G$ was obtained that came close to $i(K_{d,d})^{n/2d}$ in the second term of the exponent; this was Kahn's (unpublished) bound $i(G)\leq 2^{n/2 + O(n(1+o(1))/d)}$. This was improved to $i(G) \leq 2^{n/2 + O(n(1+o(1))/2d)}$ by Galvin \cite{Galvin-asymptot-ind-count}. Finally Zhao \cite{Zhao} deduced the exact bound for general $G$ from the bipartite case.

The aim of this note is to obtain an asymptotic version of Conjecture \ref{conj-qcols}, along the lines of Galvin's upper bound on the count of independent sets in $n$-vertex, $d$-regular graphs. Before stating the main result, we need to do some preliminary calculations. Define
$$
\eta = \eta(q) = \left\{ \begin{array}{ll} \frac{q^2}{4} & \mbox{if $q$ is even} \\ \lfloor \frac{q}{2}\rfloor \lceil \frac{q}{2} \rceil = \frac{q^2-1}{4}  & \mbox{if $q$ is odd,} \end{array}\right.
$$
and
$$
m = m(q) = \left\{ \begin{array}{ll} {q \choose q/2} & \mbox{if $q$ is even} \\  {q \choose \lfloor q/2 \rfloor} + {q \choose \lceil q/2 \rceil}  = {q+1 \choose (q+1)/2} & \mbox{if $q$ is odd.} \end{array}\right.
$$
Fix a bipartition ${\mathcal E} \cup {\mathcal O}$ of $K_{d,d}$. The set of $q$-colorings of $K_{d,d}$ may be written as $\cup_{(A,B)} {\mathcal C}(A,B)$ where the union is over ordered pairs $(A,B)$ with $|A|, |B| > 0$ and $|A\cap B|=0$, and where ${\mathcal C}(A,B)$ consists of colorings in which the set of colors appearing on ${\mathcal E}$ (resp. ${\mathcal O}$) is exactly $A$ (resp. $B$). Using inclusion-exclusion we easily get $|{\mathcal C}(A,B)| = (|A||B|)^{d/2} + O((|A||B|-1)^{d/2})$. The maximum possible value of $|A||B|$ is $\eta(q)$ (achieved when $A \cup B =\{1, \ldots q\}$ and $|A|-|B| \in \{-1,0,1\}$), and there are $m(q)$ pairs that achieve this value. It follows that
$$
\left|c_q(K_{d,d}) - \eta^d m \right| \leq (\eta - c)^d
$$
for some $c =c(q) > 0$. This leads to
\begin{equation} \label{inq-kdd}
c_q(K_{d,d})^\frac{n}{2d} = \eta^\frac{n}{2} m^\frac{n(1+o(1))}{2d}
\end{equation}
where (here and everywhere) $o(1) \rightarrow 0$ as $d \rightarrow \infty$.

Our main theorem is an upper bound on $c_q(G)$ for all $G \in {\mathcal G}(n,d)$ that matches (\ref{inq-kdd}) up to the $o(1)$ term.
\begin{thm} \label{thm-asybound}
Fix $d \geq 2$ and $n \geq d+1$. For any $G \in {\mathcal G}(n,d)$ and any $q \geq 3$,
$$
c_q(G) \leq \eta^\frac{n}{2}m^\frac{n(1+o(1))}{2d}.
$$
\end{thm}
The best previous result in this direction was from \cite{Galvin-homsconjcorr}, where it was shown that
$$
c_q(G) \leq \eta^\frac{n}{2}m^\frac{n(1-q)(1+o(1))}{dq}.
$$
(This only appears explicitly in \cite{Galvin-homsconjcorr} for $q=3$, but follows immediately for general $q$ from Proposition \ref{prop-large_ind_set} below by taking $\alpha = n/q$; note that for all smaller $\alpha$, $c_q(G)=0$.)

\medskip

To prove Theorem \ref{thm-asybound} we consider the independence number $\alpha(G)$ of $G$, the number of vertices in a largest independent set, and deal separately with large $\alpha(G)$ and small $\alpha(G)$. The case of large $\alpha(G)$ has already been dealt with in \cite[Section 5]{Galvin-homsconjcorr}, where an entropy approach is used to obtain the following result.
\begin{prop} \label{prop-large_ind_set}
Fix $d \geq 2$ and $n \geq d+1$. For any $G \in {\mathcal G}(n,d)$ with $\alpha(G) \geq n(1-\varepsilon)/2$ and any $q \geq 3$,
$$
c_q(G) \leq \eta^\frac{n}{2}m^\frac{n(1+\varepsilon)}{2d} C^\frac{n}{d^2}
$$
where $C=C(q)>1$ is a constant.
\end{prop}

To bound $c_q(G)$ when $G$ has no large independent sets we adopt an argument of Sapozhenko to obtain the following, which we prove in Section \ref{sec-smallalpha}.
\begin{lemma} \label{lem-smallindset}
Fix $d \geq 2$ and $n \geq d+1$. For each $q \geq 3$, there are constants $c_1, c_2 > 0$ such that if $G \in {\mathcal G}(n,d)$ has $\alpha(G) \leq n(1-\varepsilon)/2$ then
$$
c_q(G) \leq \eta^\frac{n}{2} \exp_2\left\{c_1n\sqrt{\frac{\log d}{d}}-c_2\varepsilon n\right\}.
$$
\end{lemma}
(For concreteness, here and throughout $\log=\log_2$.)
Along the way, we describe a very simple argument that gives the weaker bound
\begin{equation} \label{inq-weak}
c_q(G) \leq \eta^\frac{n}{2} 2^{O\left(n\sqrt{\frac{\log d}{d}}\right)}
\end{equation}
valid for {\em all} $G \in {\mathcal G}(n,d)$. (Here and elsewhere, constants implied in $O$ or $\Omega$ statements may depend on $q$.)

Taking $\varepsilon=C'\sqrt{\log d/d}$ for suitably large $C'=C'(q)>0$, Proposition \ref{prop-large_ind_set} and Lemma \ref{lem-smallindset} combine to give Theorem \ref{thm-asybound} in the following precise form: for any $G \in {\mathcal G}(n,d)$ and any $q \geq 3$,
$$
c_q(G) \leq \eta^\frac{n}{2}m^{\left(\frac{n}{2d}+O\left(\frac{n\sqrt{\log d}}{d^{3/2}}\right)\right)}.
$$
Lemma \ref{lem-smallindset} together with \cite[Prop. 1.2]{GalvinTetali-weighted} (the bipartite case of Conjecture \ref{conj-qcols}) also shows that the only $G \in {\mathcal G}(n,d)$ which remain as potential counterexamples to Conjecture \ref{conj-qcols} are those which are non-bipartite and have an independent set of size at least $n/2(1 - C\sqrt{(\log d)/d})$ for some constant $C>0$.

\medskip

A simple corollary of Lemma \ref{lem-smallindset} is that for each fixed $\varepsilon > 0$ and $q \geq 3$ there is $\delta > 0$ such that for all $d \geq 2$, $n \geq d+1$ and $G \in {\mathcal G}(n,d)$ with $\alpha(G) \leq n(1-\varepsilon)/2$, we have
$$
c_q(G) \leq (\eta-\delta)^\frac{n}{2}.
$$
A natural question to ask is how $\delta$ (more precisely, the supremum over all $\delta$ for which the preceding statement is true) varies with $\varepsilon$ in the range $0 \leq \varepsilon \leq 1-(2/q)$. At $\varepsilon =0$ we have $\delta=0$ (by Theorem \ref{thm-asybound} and the example of the disjoint union of $K_{d,d}$'s), and from the fact that $c_q(G)=0$ whenever $\alpha(G) < n/q$ we conclude that $\delta=\eta$ for all $\varepsilon > 1-(2/q)$.
\begin{question} \label{quest-threshold}
Fix $d \geq 2$, $n \geq d+1$, $q \geq 3$ and $0 \leq \varepsilon \leq 1-(2/q)$. What is the maximum of $c_q(G)$ over all $G \in {\mathcal G}(n,d)$ with $\alpha(G) \leq n(1-\varepsilon)/2$?
\end{question}

\section{Proof of Lemma \ref{lem-smallindset} --- Small independent sets} \label{sec-smallalpha}

To obtain (\ref{inq-weak}) we modify an argument due to Sapozenko \cite{Sapozhenko3}, originally used to enumerate independent sets in a regular graph; a further modification of this argument will give Lemma \ref{lem-smallindset}.

Let $\varphi = \sqrt{d \log d}/q$ (note that $\varphi < d$). For an independent set $I$ in $G$, recursively construct
sets $T(I)$ and $D(T)$ as follows. Pick $u_1 \in I$ and set $T_1 =
\{u_1\}$. Given $T_m =\{u_1, \ldots, u_m\}$, if there is $u_{m+1}
\in I$ with $N(u_{m+1})\setminus N(T_m) \geq \varphi$, then set $T_{m+1}=\{u_1, \ldots,
u_{m+1}\}$ (here $N(\cdot)$ indicates open neighborhood). If there is no such $u_{m+1}$, then set $T=T_m$ and
$$
D(T)=\{v \in V(G)\setminus N(T): N(v)\setminus N(T) < \varphi\}.
$$
Note that
\begin{equation} \label{Tbound}
|T| \leq \frac{n}{\varphi},
\end{equation}
since by construction $n \geq N(T) \geq (|T|-1)\varphi + d \geq |T|\varphi$; that
\begin{equation} \label{Dfact}
I \subseteq D
\end{equation}
since if $I \setminus D \neq \emptyset$, the
construction of $T$ would not have stopped (note that $N(T) \cap I =
\emptyset$); and that
\begin{equation} \label{Dbound}
|D| \leq \frac{nd}{2d-\varphi} \leq \frac{n}{2}\left(1+\frac{\varphi}{d}\right).
\end{equation}
The second inequality here follows from $\varphi < d$. To see the first, consider the bipartite graph with partition classes $D$
and $N(T)$ and edges induced from $G$. This graph has at most $d|N(T)| \leq d(N-|D|)$ edges
(since each vertex in $N(T)$ has at most $d$ edges to $D$, and there
are at most $N-|D|$ such vertices), and at least $(d-\varphi)|D|$ edges
(since each vertex in $D$ has at least $d-\varphi$ edges to $N(T)$).
Putting these two inequalities together gives (\ref{Dbound}).

Now a $q$-coloring of $G$ is an ordered partition of $V(G)$ into $q$ independent sets, $(I_1, \ldots, I_q)$, with $I_k$ being the set of vertices colored $k$. Following Sapozhenko's argument, we associate with this partition an ordered list $(D(T(I_1)), \ldots, D(T(I_q)))$. We recover all $q$-colorings of $G$ (and perhaps more) by finding all such lists, and then for each list $(D_1, \ldots, D_q)$ finding all ordered partitions of the $V(G)$ into $q$ sets $(I_1, \ldots, I_q)$ (not necessarily independent sets), with $I_k \subseteq D_k$ for each $k$. We say that such a partition is {\em compatible} with the $D_k$'s.

By (\ref{Tbound}) each possible $D_k$ is determined by a set of size at most $n/\varphi$ and so the number of choices for $(D_1, \ldots, D_q)$ is at most
\begin{equation} \label{int2}
\left(\sum_{i \leq n/\varphi}{n \choose i}\right)^q = 2^{O\left(n\sqrt{\frac{\log d}{d}}\right)},
\end{equation}
the equality using standard binomial estimates. We now bound the number of partitions compatible with a particular $(D_1, \ldots, D_q)$. For each $v \in V(G)$ let $a_v$ denote the number of $D_k$'s with $v \in D_k$. Using (\ref{Dbound}) we have
\begin{equation} \label{sum}
\sum_{v \in V(G)} a_v = \sum_{k=1}^q |D_k| = \frac{qn}{2}\left(1+O\left(\sqrt{\frac{\log d}{d}}\right)\right).
\end{equation}
By the AM-GM inequality we get
\begin{equation} \label{int1}
\prod_{v \in V(G)} a_v = \left(\frac{q^2}{4}\right)^\frac{n}{2} 2^{O\left(n\sqrt{\frac{\log d}{d}}\right)}.
\end{equation}
Combining (\ref{int2}) and (\ref{int1}) we get (\ref{inq-weak}) for even $q$.

We now work towards a better bound that incorporates the independence number of $G$. Since we are upper bounding $c_q$ we may assume, by adding vertices in some deterministic way if necessary, that $D_1$ satisfies $|D_1| \geq n/2$. Now we look at the subgraph induced by $D_1$. It inherits from $G$ the property that no independent set has size greater than $(n/2)(1-\varepsilon)$. This means that $D_1$ has a matching of size at least $n\varepsilon/4$ (which may be found greedily). 

Fix such a matching $M=\{x_1y_1, \ldots, x_{|M|}y_{|M|}\}$. In our naive count of colorings, we had a factor $a_{x_1}a_{x_2}$ to account for the possible colors assigned to $x_1$ and $y_1$ in a compatible partition. But since $x_1$ and $y_1$ are adjacent, we cannot assign color $1$ to both vertices, and so we have at most
$$
a_{x_1}a_{x_2}-1 = a_{x_1}a_{x_2}\left(1-\frac{1}{a_{x_1}a_{x_2}}\right) \leq a_{x_1}a_{x_2}\left(1-\frac{1}{q^2}\right)
$$
choices for this pair. Applying this argument to each of the pairs $(x_i, y_i)$, we get an upper bound on the number of colorings compatible with $(D_1, \ldots, D_q)$ of
$$
\left(\prod_{v \in V(G)} a_v\right) \left(1-\frac{1}{q^2}\right)^{|M|} = \left(\frac{q^2}{4}\right)^\frac{n}{2} 2^{O\left(n\sqrt{\frac{\log d}{d}}\right)-\Omega(\varepsilon n)},
$$
using the AM-GM inequality for the first term in the product, and our lower bound on $|M|$ for the second. Combining with (\ref{int2}) we obtain Lemma \ref{lem-smallindset} for even $q$.

Now we turn to odd $q$. Preceding exactly as before, we have
$$
c_q(G) \leq \left(\prod_{v \in V(G)} a_v\right) 2^{O\left(n\sqrt{\frac{\log d}{d}}\right)-\Omega(\varepsilon n)},
$$
so we are done (both with Lemma \ref{lem-smallindset} and with (\ref{inq-weak}) in the case of odd $q$) if we can bound
\begin{equation} \label{todo-odd}
\prod_{v \in V(G)} a_v \leq \left(\lfloor q/2 \rfloor \lceil q/2 \rceil \right)^\frac{n}{2} 2^{O\left(n\sqrt{\frac{\log d}{d}}\right)}.
\end{equation}
For this we need the following optimization lemma:
\begin{lemma} \label{lem-opt}
Let $a_1, \ldots, a_m$ be positive real numbers with average $a$. If there is a $\delta \geq 0$ such that no $a_i$ is in the interval $(a-\delta,a+\delta)$, then
$$
\prod_{i=1}^m a_i \leq \left(a^2-\delta^2\right)^\frac{m}{2} = (a-\delta)^\frac{m}{2}(a+\delta)^\frac{m}{2}.
$$
\end{lemma}

\medskip

\noindent {\em Proof}: We begin with $m$ even, say $m=2k$. Without loss of generality, assume $a_1 \leq \ldots \leq a_m$. Let ${\rm ave}_1$ be the average of $a_1$ through $a_k$ and ${\rm ave}_2$ the average of $a_{k+1}$ through $a_m$; clearly ${\rm ave}_1 \leq {\rm ave}_2$, and ${\rm ave}_1 + {\rm ave}_2 = 2a$, so ${\rm ave}_1=a-\delta'$ and ${\rm ave}_2=a+\delta'$ for some $\delta' \geq 0$. We claim that $\delta' \geq \delta$. If not, then $a_k$ (being at least ${\rm ave}_1$) is at least $a-\delta$, and so by hypothesis is at least $a+\delta$, which forces ${\rm ave}_2$ to be at least $a+\delta$, a contradiction since ${\rm ave}_2 = a + \delta' < a+\delta$.

Armed with the information that $\delta' \geq \delta$, we apply the AM-GM inequality to $a_1$ through $a_k$ and $a_{k+1}$ through $a_m$ separately and get
$$
\prod_{i=1}^m a_i  \leq (a-\delta')^k (a+\delta')^k = \left(a^2-\delta'^2\right)^\frac{m}{2} \leq \left(a^2-\delta^2\right)^\frac{m}{2}.
$$

To deal with odd $m$, we consider the problem of maximizing $\prod_{i=1}^m a_i \prod_{i=1}^m a'_i$ subject to the conditions that no $a_i$ or $a'_i$ lies in the interval $(a-\delta,a+\delta)$, and that the average of the $2m$ numbers is $a$. By the even case, the maximum is at most $(a^2-\delta^2)^m$. This remains an upper bound on the maximum if we add the conditions $a_i=a'_i$ for each $i=1, \ldots, m$; but then the maximum becomes the square of the maximum of $\prod_{i=1}^m a_i$, and we are done. \qed

\medskip

To apply Lemma \ref{lem-opt} we first assume (as we may do without loss of generality) that each $D_i$ satisfies $|D_i| \geq n/2$. This assumption together with (\ref{Dbound}) and our specific choice of $\varphi$ gives that the average of the $a_v$'s satisfies $a \in [q/2, q/2 + (1/2)\sqrt{\log d/d}]$. Since the $a_v$'s must be integers, and $\sqrt{\log d/d} < 1$, we may take $\delta = (1/2)(1-\sqrt{\log d/d})$ in Lemma \ref{lem-opt} to get
$$
\prod_{v \in V(G)} a_v \leq \left(\lfloor q/2 \rfloor + c\sqrt{\frac{\log d}{d}}\right)^\frac{n}{2}\left(\lceil q/2 \rceil \right)^\frac{n}{2} = \left(\lfloor q/2 \rfloor \lceil q/2 \rceil \right)^\frac{n}{2} 2^{O\left(n\sqrt{\frac{\log d}{d}}\right)},
$$
as required.

\end{document}